\documentclass[12pt]{amsart}
\usepackage{amssymb}
\usepackage{isolatin1}    
\usepackage[all]{xy}

\newtheorem{thm}{Theorem}[section]
\newtheorem{lem}[thm]{Lemma}
\newtheorem{cor}[thm]{Corollary}
\newtheorem{prop}[thm]{Proposition}

\newtheorem*{prob*}{Open problem}

\theoremstyle{definition}

\newtheorem{defi}[thm]{Definition}

\theoremstyle{remark}

\newtheorem{rem}[thm]{Remark}
\newtheorem*{rem*}{Remark}

\newtheorem{ex}[thm]{Example}

\DeclareMathOperator{\s}{span}

\DeclareMathOperator{\Hom}{Hom}

\newcommand{\kringel}{\mathbin{\raise1pt\hbox{$\scriptstyle\circ$}}}
\newcommand{\pkt}{\mathbin{\raise0pt\hbox{$\scriptstyle\bullet$}}}

\newcommand{\C}{\mathbb{C}}

\newcommand{\N}{\mathbb{N}}

\newcommand{\Q}{\mathbb{Q}}
\newcommand{\R}{\mathbb{R}}

\newcommand{\Der}{\mathop{\rm Der}}

\newcommand{\Lg}{\mathfrak{g}}

\newcommand{\Ll}{\mathfrak{l}}
\newcommand{\Ln}{\mathfrak{n}}

\newcommand{\Lr}{\mathfrak{r}}
\newcommand{\Ls}{\mathfrak{s}}

\newcommand{\CA}{\mathcal{A}}

\newcommand{\CH}{\mathcal{H}}
\newcommand{\CI}{\mathcal{I}}

\newcommand{\al}{\alpha}
\newcommand{\be}{\beta}
\newcommand{\ga}{\gamma}
\newcommand{\de}{\delta}

\newcommand{\la}{\lambda}
\newcommand{\om}{\omega}

\newcommand{\ra}{\rightarrow}

\renewcommand{\phi}{\varphi}

\begin{document}

\title[CNLA's and symplectic structures]{Characteristically nilpotent Lie algebras 
and symplectic structures} 

\author[D. Burde]{Dietrich Burde}
\address{Fakult\"at f\"ur Mathematik\\
Universit\"at Wien\\
  Nordbergstra\ss e 15\\
  1090 Wien \\
  Austria}
\date{\today}
\email{dietrich.burde@univie.ac.at}


\begin{abstract}
We study symplectic structures on characteristically nilpotent Lie algebras (CNLAs)
by computing the cohomology space $H^2(\Lg,k)$ for certain Lie algebras $\Lg$.
Among these Lie algebras are filiform CNLAs of dimension $n\le 14$. It turns out that there
are many examples of CNLAs which admit a symplectic structure.
A generalization of a sympletic structure is an affine structure on a Lie algebra.
\end{abstract}

\maketitle

\section{Introduction}

Invariant symplectic structures on Lie groups and on nilmanifolds
play an important role in symplectic and complex geometry.
Many questions about symplectic structures on Lie groups can be reduced
to problems in terms of the tangent Lie algebra. This leads to the study of
Lie algebras admitting a symplectic structure. Symplectic Lie algebras 
(i.e., Lie algebras admitting a symplectic structure) have been
classified in several cases. 
There is also a construction, called double extension, which yields all
symplectic nilpotent Lie algebras by successive application \cite{DAM}.
In \cite{MIL1} symplectic structures on {\it $\N$-graded} filiform Lie algebras were 
determined. 
Moreover a criterion for the existence of symplectic structures on filiform
Lie algebras was proposed.\\
In this article we study symplectic structures on {\it characteristically
nilpotent Lie algebras} (CNLAs). Such algebras do not admit an $\N$-grading since 
all derivations of CNLAs are nilpotent. Symplectic CNLAs are interesting for many
reasons. One of them is the study of Riemannian
metrics compatible with a given invariant geometric structure on a nilpotent
Lie group: in \cite{LAU} the concept of a {\it minimal} left-invariant Riemannian metric
on a nilpotent Lie group endowed with an invariant geometric structure is discussed.
For the symplectic case, such a nice minimal metric need not always exist. In fact,
there is an obstruction if the Lie algebra is symplectic and a CNLA.
Symplectic Lie algebras are also special cases of Lie algebras admitting affine structures.
Lie algebras with affine structures are the infinitesimal 
analogue of Lie groups with a left-invariant affine structure.
There have been made many efforts to solve the difficult existence question of affine
structures for a given Lie algebra \cite{BU1}, \cite{BU2}, \cite{MIL}.
From this point of view the determination of sympletic Lie algebras is also interesting.
Finally, symplectic Lie algebras play a role in superconformal field theories, see for
example \cite{PAR}. \\
The paper is organized as follows. In section two we introduce symplectic structures on
Lie algebras using Lie algebra cohomology. We recall results on the cohomology groups
$H^1$ with the coadjoint module and $H^2$ with the trivial module. We explain the relation
between affine and symplectic structures. In section three we classify all complex 
symplectic filiform CNLAs of dimension $n\le 10$. 
Here we do not use the classification of symplectic filiform Lie algebras \cite{GJKS}, 
since there are some mistakes in it. In section four we determine certain symplectic
filiform CLNAs of dimension $n\ge 12$.

\section{Preliminaries}

Let us first recall the cohomology of Lie algebras.
Let $k$ denote a field of characteristic zero and $\Lg$ a Lie algebra
over $k$.
For a $\Lg$-module $M$ the space of $p$--cochains is defined by
\begin{equation*}
C^p(\Lg,M)=\begin{cases}
\Hom_K(\Lambda^p\Lg,M) & \text{if $p\ge0$},\\
0 & \text{if $p<0$}.
\end{cases}
\end{equation*}  
The standard cochain complex $\{ C^{\pkt}(\Lg,M),d \}$
yields the space $Z^p(\Lg,M)$ of $p$-cocycles, the space
$B^p(\Lg,M)$ of $p$-coboundaries and $H^p(\Lg,M)=Z^p(\Lg,M)/B^p(\Lg,M)$,
the $p$-th cohomology space.
Let $M=k$ denote the trivial $\Lg$-module.
In that case the space of $2$--cocycles and $2$--coboundaries is given 
explicitely by
\begin{equation*}
\begin{split}
Z^2(\Lg,k) & =\{ \om \in \Hom(\Lambda^2 \Lg,k) \mid
 \om ([x_1,x_2]\wedge x_3) - \om ([x_1,x_3]\wedge x_2)\\
 & \quad + \om ([x_2,x_3]\wedge x_1)=0 \}\\
B^2(\Lg,k) & =\{ \om \in \Hom(\Lambda^2\Lg,k) \mid \om (x_1\wedge x_2)=
f([x_1,x_2])\\
 & \quad \mbox{ for some } f\in \Hom (\Lg,k)\}
\end{split}
\end{equation*}

\begin{defi}
A Lie group $G$ is said to have a left-invariant symplectic structure
if it has a left-invariant non-degenerate closed $2$-form $\om$.
\end{defi}

\begin{ex}
The Lie group $\CH_3\times \R$, where $\R$ is the abelian Lie group (with
coordinate $t$) and $\CH_3$ is the Heisenberg group
consisting of all real matrices of the form 
\[
\begin{pmatrix} 1 & x & z \\
0  & 1 & y \\
0 & 0  & 1  
\end{pmatrix}
\]
admits a left-invariant symplectic structure given by the form
\[
\om =dx\wedge (dz-xdy)+dy\wedge dt
\]
\end{ex}

\begin{defi}
A Lie algebra $\Lg$ over $k$ is called {\it symplectic}
if there is an nondegenerate $\om\in Z^2(\Lg,k)$, i.e.,  
if there exists a nondegenerate skew-symmetric bilinear form 
$B \colon \Lg \times \Lg \ra k$ satisfying
\begin{equation}\label{cocycle}
B([x,y],z)-B(x,[y,z])+B(y,[x,z])=0
\end{equation}
\end{defi}

If $\om_G$ is a left-invariant symplectic form on $G$, then
$\om_G$ defines a symplectic structure on the Lie algebra $\Lg$
of $G$. Conversely any symplectic form $\om_{\Lg}$ of $\Lg$ defines
a left-invariant symplectic structure on $G$. 
Note that a finite-dimensional symplectic Lie algebra has 
even dimension.

\begin{rem}
A symplectic Lie algebra is also called a {\it quasi-Frobenius} Lie algebra.
This is a natural generalization of a Frobenius Lie algebra. 
A Lie algebra is called {\it Frobenius} if there exists a 
nondegenerate $\om\in B^2(\Lg,k)$, i.e., a linear functional
$f\in \Hom(\Lg,k)$ such that $B(x,y)=f([x,y])$ is nondegenerate.
Hence Frobenius Lie algebras are symplectic. They
have been studied in various contexts, see \cite{ELA}, \cite{OOM}. 
Many properties are known: they have trivial center, no non-zero semisimple ideals
and a non-nilpotent solvable radical, see \cite{OOM}.
Moreover, a Lie algebra $\Lg$ of a linear algebraic group $G$
over an algebraically closed field of characteristic zero 
is Frobenius if and only if the universal enveloping algebra $U(\Lg)$ is primitive, 
and if and only if $G$ admits an open orbit in the coadjoint module. \\
Quasi-Frobenius Lie algebras have been studied
in connection with rational solutions of the classical Yang-Baxter
equation (CYBE) \cite{STO} (there is a correspondence between
rational solutions of CYBE for a simple Lie algebra $\Lg$ and 
quasi-Frobenius subalgebras of $\Lg$). They appear also
in superconformal field theories (see \cite{PAR}) and related subjects.
\end{rem}

\begin{ex}
Clearly any abelian Lie algebra of even dimension is symplectic.\\
In dimension $2$ over the complex numbers there are two Lie algebras,
$\C^2$ and the non-abelian Lie algebra $\Lr_2(\C)$, which is given
by $[e_1,e_2]=e_1$ where $(e_1,e_2)$ denotes a basis of $\C^2$.
The algebra $\Lr_2(\C)$ is Frobenius, and hence symplectic.
\end{ex}

\begin{ex}\label{n4}
Let $\Ln_4$ be the $4$-dimensional nilpotent Lie algebra with
basis $(e_1,e_2,e_3,e_4)$ defined by the brackets
\begin{equation*}
[e_1,e_2]=e_3, \; [e_1,e_3]=e_4
\end{equation*}
Clearly this Lie algebra is not Frobenius since it has a non-trivial center.
It is easy to see that the space $H^2(\Ln_4,k)$ is spanned by the classes of
$\om_1$ and $\om_2$ which are defined by
\begin{align*}
\om_1(e_1\wedge e_4)& = 1 \\
\om_2(e_2\wedge e_3)& = 1
\end{align*} 
With respect to the given basis, the matrix of $\om_1+\om_2$ 
is given by
\begin{equation*}
\begin{pmatrix} 0 & 0 & 0 & 1\\
0  & 0  & 1 & 0 \\
0  & -1 & 0 & 0 \\
-1 & 0  & 0 & 0 
\end{pmatrix}
\end{equation*}     
Since $\om_1+\om_2$ is nondegenerate, $\Ln_4$ is quasi-Frobenius, or symplectic.
\end{ex}

There is a large literature on symplectic Lie algebras. In \cite{DAM} it was shown
that all symplectic nilpotent Lie algebras can be obtained by a consecutive
procedure, called double extensions, starting with the Lie algebra $\{0\}$.
The classification of complex symplectic filiform Lie algebras of dimension $n\le 10$
up to symplecto-isomorphism was given in \cite{GJKS}. However, there are some mistakes
in it.
Also, the classification of all symplectic Lie algebras in dimension $n\le 4$ is well known. 
Let us recall the result for $n=4$ over the complex numbers.

\begin{prop}
Any $4$-dimensional complex quasi-Frobenius Lie algebra is 
isomorphic to one and only one Lie algebra of the following list:

\vspace*{0.5cm}
\begin{center}
\begin{tabular}{c|c}
 $\Lg$ & Defining Lie brackets \\
\hline     
$\C^4$ & $-$ \\ 
$\Ln_3(\C)\oplus \C$ & $[e_1,e_2]=e_3$ \\ 
$\Lr_2(\C) \oplus  \C^2$ & $[e_1,e_2]=e_1$ \\ 
$\Lr_{3,-1}(\C) \oplus \C $ & $[e_1,e_2]=e_2, [e_1,e_3]=- e_3$ \\ 
$\Lr_2(\C) \oplus \Lr_2(\C)$ & $[e_1,e_2]=e_1, [e_3,e_4]=e_3$ \\  
$\Ln_4(\C)$ & $[e_1,e_2]=e_3, [e_1,e_3]=e_4$ \\ 
$\Lg_1(-1)$ & $ [e_1,e_2]=e_2, [e_1,e_3]=e_3, [e_1,e_4]=-e_4$ \\  
$\Lg_2(\al,\al)$ & $ [e_1,e_2]=e_3, [e_1,e_3]=e_4, [e_1,e_4]=\al e_2 -
\al e_3 +e_4$ \\  
$\Lg_6$ & $[e_1,e_2]=e_2, [e_1,e_3]=e_3, [e_1,e_4]=2e_4, [e_2,e_3]=e_4$\\  
$\Lg_8(\al)$ & $[e_1,e_2]=e_3, [e_1,e_3]=-\al e_2+e_3,[e_1,e_4]=e_4,
[e_2,e_3]=e_4$ \\   
\end{tabular}
\end{center}
\end{prop}

\vspace*{0.5cm}
\begin{proof}
Using the classification of $4$-dimensional Lie algebras given in
\cite{BS} we determine the spaces $Z^2(\Lg,\C)$. The result follows 
by computing the determinants. We want to demonstrate the details by
taking one example, the Lie algebra $\Lg=\Lr_{3,\la}\oplus\C$. It is defined by
the brackets $[e_1,e_2]=e_2, [e_1,e_3]=\la e_3$ where $\la \in \C, 0<|\la|\le
1$. The space $Z^2(\Lg,\C)$ is represented by the subspace
of matrices of the form
\begin{equation*}
\begin{pmatrix} 0 & \al & \be & \ga\\
-\al  & 0  & \de & 0 \\
-\be  & -\de & 0 & 0 \\
-\ga & 0  & 0 & 0 
\end{pmatrix}
\end{equation*} 
with determinant $(\ga \de)^2$ and the condition
$(\la +1)\de =0$. Hence $\Lg$ is symplectic if and only if $\la=-1$.
Note that the Lie algebra $\Ls\Ll_2(\C)\oplus \C$ is not symplectic.
It follows that all $4$-dimensional symplectic Lie algebras are solvable. 
\end{proof}

A further generalization of a symplectic Lie algebra is a Lie algebra
admitting an affine structure.

\begin{defi}
A vector space $A$ over $k$ together with a $k$-bilinear product
$A \times A \rightarrow A, (x,y) \mapsto x\cdot y$ is called 
{\it left-symmetric algebra} or LSA, if
\begin{equation}\label{lsa1}
x\cdot (y\cdot z)-(x\cdot y)\cdot z= y\cdot (x\cdot z)-(y\cdot x)\cdot z
\end{equation}
for all $x,y,z \in A$. 
\end{defi}

The left-multiplication $L$ in $A$ is given by $L(x)y=x\cdot y$.

\begin{defi}\label{affine}
An {\it affine structure} on a Lie algebra
$\Lg$ over $k$ is a $k$--bilinear product $\Lg \times \Lg \rightarrow \Lg$
satisfying $(\ref{lsa1})$ and 
\begin{equation}\label{lsa2}
[x,y]=x\cdot y -y\cdot x
\end{equation}
for all $x,y,z \in \Lg$. A Lie algebra over $k$ admitting an affine
structure is also called {\it affine}.
\end{defi}       

The term affine Lie algebra is also used differently in the literature.
Note that we have two different bilinear products in the above definition:
the Lie bracket and the dot product. \\
The conditions may be reformulated as follows: 
there exists a bilinear product $x\cdot y$ on $\Lg \times \Lg$ which 
defines a $\Lg$-module structure on $\Lg$ itself, denoted by $\Lg_L$, such 
that the identity mapping $\iota: \Lg \ra \Lg_L$ is a 1-cocycle in 
$Z^1(\Lg,\Lg_L)$. In other words, $(\ref{lsa1})$
and $(\ref{lsa2})$ are equivalent to the following identities:
\begin{align*}
[x,y]\cdot z & = x\cdot(y\cdot z)-y\cdot(x\cdot z)\\
\iota([x,y]) & = x\cdot \iota(y)-y\cdot \iota(x)
\end{align*}

In general it is very difficult for a given Lie algebra to decide
whether it is affine or not. It is well known that a Lie
algebra $\Lg$ over charactristic zero satisfying $\Lg=[\Lg,\Lg]$
is not affine. Hence the existence problem mainly arises for
solvable Lie algebras. In low dimensions this problem has a positive
solution. All complex Lie algebras of dimension $n\le 4$ are affine except
for $\Ls\Ll_2(\C)$, and all complex nilpotent Lie algebras of dimension
$n\le 7$ are affine.
There exist already examples of nilpotent Lie algebras of dimension $10$
which are not affine, see \cite{BU2}. Geometrically this means that
there are nilmanifolds which are not affine.

\begin{lem}\label{inv}
A Lie algebra $\Lg$ is affine if and 
only if there exists a $\Lg$-module structure $M$ on the vector
space of $\Lg$ such that there is a linear map 
$\phi \in Z^1(\Lg,M)$ satisfying $\det \phi \ne 0$.
\end{lem}

\begin{proof}
Let $\phi$ be a nonsingular $1$-cocycle in $Z^1(\Lg,M)$ and denote
the action of $\Lg$ on $M$ by $(x,m)\mapsto x \pkt m$. Then define
a bilinear product on $\Lg$ by
\begin{equation*}
x\cdot y:=\phi^{-1}(x \pkt \phi (y))
\end{equation*}
This product is left-symmetric, since it defines a $\Lg$-module
structure on $\Lg$, obtained by conjugation with $\phi$ from $M$,
satisfying $x\cdot y-y\cdot x=\phi^{-1}(x \pkt \phi (y)-y \pkt \phi (x))
= \phi^{-1}(\phi([x,y]))=[x,y]$.
Conversely, an affine structure on $\Lg$ yields a nonsingular
$1$-cocycle $\iota \in Z^1(\Lg,\Lg_L)$.
\end{proof}

Denote by $\Lg$ the adjoint module and by $\Lg^*=\Hom(\Lg,k)$ the
coadjoint module of $\Lg$. The coadjoint action is given by
$(x,f)\mapsto x\pkt f$ where $(x\pkt f)y=-f([x,y])$ for $x\in \Lg$ and
$f \in \Lg^*$. A derivation $D\in \Der (\Lg)$ is just a $1$--cocycle in 
$Z^1(\Lg,\Lg)$. The following two corollaries are easily derived from
lemma $\ref{inv}$: 

\begin{cor}
Any Lie algebra $\Lg$ admitting a nonsingular $D\in \Der (\Lg)$ is affine. 
\end{cor}

\begin{cor}\label{dual}
Any Lie algebra $\Lg$ admitting a nonsingular $\phi\in Z^1(\Lg,\Lg^*)$
is affine.
\end{cor}

For the cohomology with coefficients in the dual module we have
the following well-known result.

\begin{prop}\label{sub}
$H^2(\Lg,k)$ may be regarded as a subspace of $H^1(\Lg,\Lg^*)$. 
If $\Lg$ does not have a non-zero invariant bilinear form then
$H^2(\Lg,k)\simeq H^1(\Lg,\Lg^*)$.
\end{prop}

\begin{proof}
The space $Z^1(\Lg,\Lg^*)$ may be interpreted as the space of bilinear forms
$B: \Lg \times \Lg \ra k$ satisfying condition \eqref{cocycle}, i.e., 
$B([x,y],z)-B(x,[y,z])+B(y,[x,z])=0$. Indeed, if $\phi \in Z^1(\Lg,\Lg^*)$, then
define $B$ by $B(x,y)=\phi(x)y$. The condition 
$\phi([x,y])=x\pkt \phi(y)-y\pkt \phi(x)$ is just equivalent to the 
condition \eqref{cocycle} on $B$. Conversely, given a $B$ which sastifies 
identity \eqref{cocycle}, $\phi$ defined by $\phi(x)y=B(x,y)$ will be a 
$1$-cocycle in $Z^1(\Lg,\Lg^*)$.
Now the subspace formed by those $B$ which are skew-symmetric, i.e., satisfy
$B(x,y)=-B(y,x)$, corresponds exactly to the space $Z^2(\Lg,k)$.
Since obviously $B^2(\Lg,k) \simeq B^1(\Lg,\Lg^*)$, $H^2(\Lg,k)$ becomes
a subspace of $H^1(\Lg,\Lg^*)$. 
To prove the second claim, let $B$ be a bilinear form in 
$Z^1(\Lg,\Lg^*)$ and define $B^*$ by $B^*(x,y)=B(y,x)$. Then $\be=B+B^*$ is an 
invariant bilinear form on $\Lg$, i.e., satisfies
$\be([x,y],z)=\be(x,[y,z])$. 
Assume that $\Lg$ does not have a non-zero invariant bilinear form. 
Then $\be$ is zero, hence $B$ is skew-symmetric and contained in
$Z^2(\Lg,k)$. 
\end{proof}

We obtain the following corollary.

\begin{cor}
Any symplectic Lie algebra is affine.
\end{cor}

\begin{proof}
Let $\om \in Z^2(\Lg,k)$ be nondegenerate. Then by Proposition $\ref{sub}$, $\phi$ 
defined by $\phi(x)y=\om (x\wedge y)$ is a nonsingular $1$-cocycle in
$Z^1(\Lg,\Lg^*)$ since $\ker (\phi)=\{x\in \Lg \mid \om(x\wedge y)=
0 \text{ for all } y\in \Lg\}=0$. The claim follows from Corollary $\ref{dual}$.
\end{proof}

The corollary is well known. A different proof can be found in \cite{GJKS}.
Clearly an affine Lie algebra need not be symplectic, since there
exist affine Lie algebras of odd dimension. There are also affine
Lie algebras of even dimension which are not symplectic (see the
example below). Although Lie algebras of odd dimension admit no nonsingular
$\om \in Z^2(\Lg,k)$, there may be a nonsingular $\phi\in Z^1(\Lg,\Lg^*)$.
Easy examples are the Heisenberg Lie algebra or a filiform Lie algebra 
of dimension $5$. Hence there exist nilpotent Lie algebras with a non-zero
invariant bilinear form. 

\begin{ex}\label{contre}
Let $\Lg$ be the $6$-dimensional nilpotent Lie algebra with
basis $(e_1,\ldots, e_6)$ defined by the brackets
\begin{align*}
[e_1,e_i] & = e_{i+1}, \; 2\le i\le 5 \\
[e_2,e_5] & = -e_6 \\
[e_3,e_4] & = e_6
\end{align*}
Then $H^2(\Lg,k)$ is spanned by the classes of $\om_1$ and $\om_2$,
which are defined by 
\begin{align*}
\om_1(e_2\wedge e_3)& = 1 \\
\om_2(e_2\wedge e_5)& = 1,\;\om_2(e_3\wedge e_4) = -1
\end{align*}  
For $\om \in Z^2(\Lg,k)$ define $\phi_{\om}$ by $\phi_{\om}(x)y=
\om(x\wedge y)$. If $z\in Z(\Lg)$, then $\phi_{\om}(z)=0$ for any
$\om \in B^2(\Lg,k)$. In our case $Z(\Lg)$ is spanned by $e_6$, and
any linear combination $\om$ of $\om_1$ and $\om_2$ satisfies
$\phi_{\om}(e_6)=0$. Hence there is no nondegenerate $\om \in Z^2(\Lg,k)$
and the Lie algebra is not symplectic. Nevertheless $\Lg$ admits an
affine structure induced by a nonsingular derivation. In fact, 
it is easy to see that the linear map $d : \Lg \ra \Lg$
given by $d(e_i)=ie_i, i=1,\ldots,5$ and $d(e_6)=7e_6$ defines a 
nonsingular derivation. 
\end{ex}

\section{Symplectic filiform CNLAs of dimension $n<12$}

In this section we classify characteristically nilpotent
symplectic filiform Lie algebras of dimension $n< 12$. 
There is a classification of complex symplectic filiform Lie algebras of
dimension $n<12$ up to symplecto-isomorphism \cite{GJKS}.
However, there are some mistakes in it; see also the remark in \cite{MIL}.
We use a different method which does not rely on the explicit classification:
in \cite{BU1} we have computed the cohomology space $H^2(\Lg,k)$
for filiform nilpotent Lie algebras. Consequently we can use
the knowledge of $Z^2(\Lg,k)$ to determine symplectic filiform Lie algebras.

\begin{defi}
Let $\Lg$ be a nilpotent Lie algebra and $\{ \Lg^k\}$ its lower central series
defined by $\Lg^0=\Lg,\;\Lg^k=[\Lg^{k-1},\Lg]$ for $k\ge 1$.
There exists an integer $p$ such that $\Lg^p=0$ and $\Lg^{p-1}\neq 0$,
called nilindex of $\Lg$. A nilpotent Lie algebra of dimension $n$
and nilindex $p=n-1$ is called {\it filiform}.
\end{defi} 

We divide the set $\CA_n$ of filiform Lie algebra laws of dimension
$n$ into subsets $\CA_{n,i}$ such that algebras from different subsets
are non-isomorphic (but may be isomorphic if they belong to the same
subset), and algebras belonging to the same subset have the same
second scalar cohomology. 
If $\Lg$ is a filiform Lie algebra of dimension $n$,
then there exists an adapted basis $(e_1,\ldots,e_n)$ for $\Lg$, see
\cite{BU1}. We write $\CA_n$ for the set of elements which are the structure
constants of a filiform Lie algebra with respect to an adapted basis.
The brackets of such a filiform Lie algebra
with respect to the basis $(e_1,\ldots,e_n)$ are then given by

\begin{align}\label{lie}
[e_1,e_i] & =e_{i+1}, \quad i=2,\dots ,n-1 \\
[e_i,e_j] & =\sum_{r=1}^n\biggl(\;\sum_{\ell=0}^{[(j-i-1)/2]} (-1)^\ell
{ j-i-\ell-1 \choose \ell}\al_{i+\ell,\, r-j+i+2\ell+1}\biggr)e_r,
 \quad 2 \le i<j \le n.
\end{align}
with constants $\al_{k,s}$ which are zero for all pairs $(k,s)$ not in 
the index set $\CI_n$. Here 
$\CI_n$ is given by
\begin{align*}
\CI_n^0 &=\{(k,s)\in \N \times \N \mid 2 \le k \le [n/2],\,
2k+1 \le s \le n \},\\
\CI_n& =\begin{cases}
\CI_n^0 & \text{if $n$ is odd},\\
\CI_n^0 \cup \{(\frac{n}{2},n)\} & \text{if $n$ is even}.
\end{cases}
\end{align*}     

Let $f = 3\al_{4,10}(\al_{2,6}+\al_{3,8})-4\al_{3,8}^2$.
The above mentioned subsets $\CA_{n,i}$ are given as follows, for $n=4,6,8,10$:

\vspace*{0.5cm}
\begin{center}
\begin{tabular}{c|c}
Class & Conditions \\
\hline
$\CA_{4,1}$ & $-$ \\ 
$\CA_{6,1}$ & $\al_{3,6}\ne 0$ \\ 
$\CA_{6,2}$ & $\al_{3,6}= 0$   \\ 
$\CA_{8,1}$ & $\al_{4,8}\ne 0,\,2\al_{2,5}+\al_{3,7}=0$   \\ 
$\CA_{8,2}$ & $\al_{4,8}=0,\,2\al_{2,5}+\al_{3,7}\ne 0$   \\ 
$\CA_{8,3}$ & $\al_{4,8}=0,\,2\al_{2,5}+\al_{3,7}=0,\, \al_{2,5}\ne 0$   \\ 
$\CA_{8,4}$ & $\al_{2,5}=\al_{3,7}=\al_{4,8}=0$   \\  
$\CA_{10,1}$ & $\al_{5,10}\ne 0,\,2\al_{2,5}+\al_{3,7}\ne 0$ \\ 
$\CA_{10,2}$ & $\al_{5,10}\ne 0,\,2\al_{2,5}+\al_{3,7}=   0$ \\ 
$\CA_{10,3}$ & $\al_{5,10}= 0,\,2\al_{2,5}+\al_{3,7}\ne 0,\, \al_{3,7}^2\ne \al_{2,5}^2 $ \\ 
$\CA_{10,4}$ & $\al_{5,10}= 0,\,2\al_{2,5}+\al_{3,7}\ne 0,\, \al_{3,7}^2= \al_{2,5}^2$ \\ 
$\CA_{10,5}$ & $\al_{5,10}= 0,\,2\al_{2,5}+\al_{3,7}= 0,\, \al_{4,9}\ne 0,\, \al_{2,6}^2+2\al_{2,7}
\al_{4,9}\ne 0$ \\ 
$\CA_{10,6}$ & $\al_{5,10}= 0,\,2\al_{2,5}+\al_{3,7}= 0,\, \al_{4,9}\ne 0,\,\al_{2,6}^2+2\al_{2,7}
\al_{4,9}= 0$ \\ 
$\CA_{10,7}$ & $\al_{5,10}= 0,\,2\al_{2,5}+\al_{3,7}= 0,\,\al_{4,9}=0,\, 2\al_{2,7}+\al_{3,9}\ne 0$ \\ 
$\CA_{10,8}$ & $\al_{5,10}= 0,\,2\al_{2,5}+\al_{3,7}= 0,\,\al_{4,9}=0,\, 
2\al_{2,7}+\al_{3,9}= 0,\, f\ne 0$ \\ 
$\CA_{10,9}$ & $\al_{5,10}= 0,\,2\al_{2,5}+\al_{3,7}=
0,\,\al_{4,9}=0,\,2\al_{2,7}+
\al_{3,9}= 0,\, f= 0$ \\    
\end{tabular}
\end{center}

\vspace*{0.5cm}   

\begin{defi}
Let $\Lg$ be a filiform Lie algebra and $(e_1,\ldots,e_n)$ be an adapted
basis of $\Lg$. Define 
$\om_{\ell}\in \Hom (\Lambda^2\Lg,k)$ by
\begin{align}\label{omega}
\om_{\ell}(e_k\wedge e_{2\ell+3-k}) & = (-1)^k \quad \text{ for } 1\le \ell
\le [(n-1)/2],\; 2\le k\le \left[(2\ell+3)/2\right]
\end{align}  
where the values not defined (and which are not a consequence of skew-symmetry)
are understood to be zero.
\end{defi}

In general, the $\om_{\ell}$ need not be cocycles for $\ell \ge 3$.
On the other hand we know the following \cite{BU1}:

\begin{lem}
Let $\Lg$ be filiform of dimension $n\ge 5$. Then 
$\om_1,\om_2\in Z^2(\Lg,k)$. Any $2$-coboundary $\be \in B^2(\Lg,k)$
is degenerate. If $\ell<[(n-1)/2]$, then 
$\om_{\ell}$ is degenerate.
\end{lem}

In fact, for all $x\in \Lg$ and $z \in Z(\Lg)$ we have 
$\be(x\wedge z)=f([x,z])=0$ for some linear form $f\in \Hom(\Lg,k)$. 
Recall that the center $Z(\Lg)$ is $1$-dimensional.
Likewise $\om_{\ell}$ is zero on $\Lg\wedge Z(\Lg)$ for
$\ell<[(n-1)/2]$.\\[0.5cm]
Let $n=4$: the cohomology does not depend on the structure constants. 
Over the complex numbers there is only one filiform Lie algebra, 
namely $\Lg=\Ln_4(\C)$. It is not a CNLA. The result is as follows, see
example \eqref{n4}:

\begin{prop}
We have $$H^2(\Ln_4,\C)=\s \{[\om_1],[\om]\}$$
where the $2$-cocycles are defined by 
$\om_1(e_2\wedge e_3) = 1 $ and $\om(e_1\wedge e_4) = 1 $.
Since $\om +\om_1$ is nondegenerate, $\Ln_4$ is symplectic. 
\end{prop}

Let $n=6$. Denote by $\la\in \CA_6$ the law of $\Lg$. It is well known
that all such $\la$ are $\N$-graded. Hence they are not CNLAs.

\begin{prop}
In dimension $6$ we have
\begin{equation*}
H^2(\Lg,k)=\begin{cases}
\s\{[\om_1],[\om_2]\} & \text{if $\la\in\CA_{6,1}$}\\[0.2cm]
\s\{[\om_1],[\om_2],[\om]\} & \text{if $\la\in\CA_{6,2}$}\\[0.2cm]
\end{cases}
\end{equation*}
If  $\la\in\CA_{6,1}$ then $\Lg$ is not symplectic. If $\la\in\CA_{6,2}$
then $\Lg$ is symplectic.
\end{prop}

\begin{proof}
The $2$-cocycles $\om_1,\om_2$ are defined as in \eqref{omega}, and $\om$
is defined by

\begin{align*}
\om(e_1\wedge e_6)& = 1 \\
\om(e_3\wedge e_4)& = \al_{2,5} \\
\om(e_2\wedge e_4)& = \al_{2,6} \\
\end{align*}

Computing the determinant we obtain
$$\det (r \om_2+s \om)=(r-s \al_{2,5})^2r^2s^2$$ 
which is non-zero for a suitable choice of the constants $r$ and $s$.
Hence all $\la\in \CA_{6,2}$ are symplectic.
\end{proof}

Using the classification list and the notation of \cite{GJK} over $\C$ we obtain:

\begin{cor}
Every $6$-dimensional complex filiform symplectic Lie algebra is isomorphic
to one of the following:

\begin{align*}
\mu_6^1 & \colon \mu_0 \\
\mu_6^2 & \colon \mu_0+ \psi_{2,5} \\
\mu_6^3 & \colon \mu_0+ \psi_{2,6} 
\end{align*}
\end{cor}

Note that our  $\psi_{i,j}$ correspond to $\Psi_{i-1,j-1}$ in \cite{GJK}.
For example, the brackets of $\mu_6^2$ are given by \eqref{lie}, $(5)$
with $\al_{2,5}=1, \al_{2,6}=0, \al_{3,6}=0$.\\[0.5cm]

Let $n=8$. Denote by $\la\in \CA_8$ the law of $\Lg$.

\begin{prop}
We have
\begin{equation*}
H^2(\Lg,k)=\begin{cases}
\s\{[\om_1],[\om_2],[\om_3]\} & \text{if $\la\in\CA_{8,1}$ or 
$\la\in\CA_{8,3}$}\\[0.2cm]
\s\{[\om_1],[\om_2],[\om]\} & \text{if $\la\in\CA_{8,2}$}\\[0.2cm]
\s\{[\om_1],[\om_2],[\om_3],[\be]\} & \text{if $\la\in\CA_{8,4}$}\\[0.2cm]
\end{cases}
\end{equation*}
where $\om$ is defined by

\begin{align*}
\om(e_1\wedge e_8)& = 1, \\
\om(e_2\wedge e_4)& = \al_{2,8},\quad \om(e_2\wedge e_6) = \al_{2,6}-2\al_{3,8},\quad
\om(e_2\wedge e_7)=\frac{\al_{2,5}(2\al_{2,5}-5\al_{3,7})}{2\al_{2,5}+
  \al_{3,7}},\\
\om(e_3\wedge e_4)& =  \al_{3,7},\quad \om(e_3\wedge e_5) =\al_{3,8},\quad
\om(e_3\wedge e_6) =\frac{2\al_{3,7}(\al_{2,5}-\al_{3,7})}{2\al_{2,5}+ \al_{3,7}}, \\
\om(e_4\wedge e_5)& = \frac{3\al_{3,7}^2}{2\al_{2,5}+ \al_{3,7}} \\
\end{align*}

and $\be$ is defined by
\begin{align*}
\be(e_1\wedge e_8)& = 1 \\
\be(e_2\wedge e_4)& = \al_{2,8},\;\be(e_2\wedge e_6) = \al_{2,6}-2\al_{3,8},\;
\be(e_2\wedge e_7)= 1\\
\be(e_3\wedge e_4)& =  \al_{3,7},\;\be(e_3\wedge e_5) =\al_{3,8},\;
\be(e_3\wedge e_6) =-1 \\
\be(e_4\wedge e_5)& = 1 \\
\end{align*}
\end{prop}

Note that $\be$ is non-degenerate. Computing determinants we obtain:

\begin{cor}
If $\la\in\CA_{8,1}$ or $\CA_{8,3}$, then $\Lg$
is not symplectic. If $\la\in\CA_{8,4}$, then $\Lg$ is symplectic.
If $\la\in\CA_{8,2}$, then $\Lg$ is symplectic if and only if
$$\al_{2,5}\al_{3,7}(\al_{2,5}-\al_{3,7})(5\al_{3,7}-2\al_{2,5}) \ne 0$$
\end{cor}

Again using the classification list of \cite{GJK} we obtain:

\begin{cor}
Every $8$-dimensional complex symplectic filiform Lie algebra is isomorphic
to one of the following laws:

\begin{align*}
\mu_8^5(\al) & \colon \mu_0 +\al \psi_{2,5}+\psi_{3,7}+\psi_{3,8},\quad \al \neq \textstyle{-\frac{1}{2},0,
1,\frac{5}{2}} \\
\mu_8^6(\al) & \colon \mu_0 +\al \psi_{2,5}+\psi_{3,7},\quad \al \neq \textstyle{-\frac{1}{2},0,1,
\frac{5}{2}} \\
\mu_8^9(\al) & \colon \mu_0 +\al \psi_{2,6}+\psi_{2,7}+\psi_{3,8} \\
\mu_8^{10}(\al) & \colon \mu_0 +\al \psi_{2,6}+\psi_{3,8} \\
\mu_8^{11}(0) & \colon \mu_0 +\psi_{2,7}+\psi_{2,8} \\
\mu_8^{15} & \colon \mu_0 +\psi_{2,6}+\psi_{2,7} \\
\mu_8^{16} & \colon \mu_0 +\psi_{2,6} \\
\mu_8^{17} & \colon \mu_0 +\psi_{2,7} \\
\mu_8^{18} & \colon \mu_0 +\psi_{2,8} \\
\mu_8^{19} & \colon \mu_0
\end{align*}
\end{cor}

It is not difficult to compute the derivations of these algebras. This yields

\begin{cor}
Every $8$-dimensional complex symplectic filiform CNLA is isomorphic to one of
the following laws: $\mu_8^5(\al), \al \neq -\frac{1}{2},0,1,\frac{5}{2}$, or $\mu_8^9(\al),
\mu_8^{11}(0), \mu_8^{15} $.
\end{cor}

Let $n=10$. Denote by $\la\in \CA_{10}$ the law of $\Lg$. 

\begin{prop}
We have

\begin{equation*}
H^2(\Lg,k)=\begin{cases}
\s\{[\om_1],[\om_2],[\om_3]\} & \text{if $\la\in\CA_{10,1}$, 
$\CA_{10,4}$ or $\CA_{10,5}$}\\[0.2cm]
\s\{[\om_1],[\om_2],[\om_3],[\om_4]\} & \text{if $\la\in\CA_{10,2}$ or 
$\CA_{10,8}$}\\[0.2cm]
\s\{[\om_1],[\om_2],[\om]\} & \text{if $\la\in\CA_{10,3}$}\\[0.2cm]
\s\{[\om_1],[\om_2],[\om_3],[\be_1]\} & \text{if $\la\in\CA_{10,6}$}\\[0.2cm]
\s\{[\om_1],[\om_2],[\om_3],[\be_2]\} & \text{if $\la\in\CA_{10,7}$}\\[0.2cm]
\s\{[\om_1],[\om_2],[\om_3],[\om_4],[\be_3]\} & \text{if $\la\in\CA_{10,9}$}\\[0.2cm]
\end{cases}
\end{equation*}
\end{prop}

The cocycles $\om,\be_1,\be_2,\be_3$ are too complicated to be listed here.
Let $p(x,y)$ denote the following polynomial

\begin{equation*}
\begin{split}  
p(x,y) & = (5y^3 - 8y^2x + 16yx^2 - 4x^3)
 (5y^3 - 16y^2x + 10yx^2 - 2x^3)\\
 & \quad \; (5y^2 - 4yx + 2x^2)(7y - 4x)y
\end{split}
\end{equation*} 
A straightforward computation of determinants yields the following result:

\begin{cor}
A filiform Lie algebra with law in
$\CA_{10,1},\CA_{10,2},\CA_{10,4},\CA_{10,5},\CA_{10,6},\CA_{10,8}$ 
is not symplectic.
Any filiform Lie algebra with law in $\CA_{10,9}$ is symplectic.
An algebra with law in $\CA_{10,3}$ is symplectic if and only if 
$p(\al_{2,5},\al_{3,7})\neq 0$.
An algebra with law in $\CA_{10,7}$ is symplectic if and only
if $f=3\al_{4,10}(\al_{2,6}+\al_{3,8})-4\al_{3,8}^2 \neq 0$.
\end{cor}

Using the classification list of \cite{GJK} we obtain:

\begin{cor}
Every $10$-dimensional complex symplectic filiform Lie algebra is isomorphic
to one of the following laws: $\mu_{10}^9(\al), \;  \mu_{10}^{10},\;  \mu_{10}^{11},\;\mu_{10}^{13}(\al,\be),\; 
\mu_{10}^{16}(\al,\be),\; \mu_{10}^{26}(0,\be),\;  \mu_{10}^{27},\;
\mu_{10}^{29}(\al),$ $\mu_{10}^{32}(-\frac{1}{2},\be),\; \mu_{10}^{33}(-\frac{1}{2}),\;
\mu_{10}^{34}(-\frac{1}{2}),\; \mu_{10}^{36}(0,\be),\;\mu_{10}^{37}(\al),\;\mu_{10}^{39}(\al),\;
\mu_{10}^{40},\;\mu_{10}^{41},\;\mu_{10}^{45}(\al),\;\mu_{10}^{46},\; \mu_{10}^{47},\; \mu_{10}^{48},\; 
 \mu_{10}^{49},$ $ \mu_{10}^{50},\; \mu_{10}^{51}$ or

\begin{align*}
\mu_{10}^1(\al,\be) & \colon \quad (\al+2)(\al+1)(\al-1)p(1,\al) \neq 0 \\
\mu_{10}^2(\al,\be) & \colon \quad (\al+2)(\al+1)(\al-1)p(1,\al) \neq 0 \\
\mu_{10}^3(\al) & \colon \quad (\al+2)(\al+1)(\al-1)p(1,\al) \neq 0 \\
\mu_{10}^4(\al) & \colon \quad (\al+2)(\al+1)(\al-1)p(1,\al) \neq 0 \\
\mu_{10}^{12}(\al,\be) & \colon \quad \al\neq 0 \\ 
\mu_{10}^{17}(\al,\be,\ga) & \colon \quad (2\be+1)(3\al+3\ga-4\ga^2) \neq 0 \\ 
\mu_{10}^{17}(\al,\be,\ga) & \colon \quad 2\be+1=3\al+3\ga-4\ga^2= 0 \\
\mu_{10}^{18}(\al,\be,\ga,\de) & \colon \quad (2\be+\de)(3\al+3\ga-4\ga^2)\neq 0 \\ 
\mu_{10}^{18}(\al,\be,\ga,\de) & \colon \quad 2\be+\de=3\al+3\ga-4\ga^2=\ga(\ga-6)= 0 \\ 
\mu_{10}^{19}(\al,\be) & \colon \quad 3\al+3\be-4\be^2 \neq 0 \\
\mu_{10}^{20}(\al,\be,\ga) & \colon \quad (2\be+\ga)(6\al+1)\neq 0 \\
\mu_{10}^{20}(\al,\be,\ga) & \colon \quad 2\be+\ga=6\al+1= 0 \\
\mu_{10}^{21}(\al,\be,\ga) & \colon \quad \ga(2\al+\be)\neq 0 \\
\mu_{10}^{22}(\al,\be,\ga) & \colon \quad \ga(2\al+\be)\neq 0 \\
\mu_{10}^{23}(\al,\be) & \colon \quad 3\al+3\be-4\be^2 = 0 \\
\mu_{10}^{24}(\al,\be,\ga) & \colon \quad 2\al+\be= 0, \al\neq 0 \\
\mu_{10}^{25}(\al,\be) & \colon \quad \al\neq 0 \\
\mu_{10}^{28}(\al,\be) & \colon \quad 2\be +1\neq 0 \\
\mu_{10}^{31}(\al,\be,\ga) & \colon \quad 2\al +\be= 0 \\
\end{align*}
\end{cor}

\begin{cor}
Every $10$-dimensional complex symplectic filiform CNLA is isomorphic
to one of the following laws: $\mu_{10}^9(\al), \; \mu_{10}^{10},\; \mu_{10}^{13}(\al,\be),\;   
\mu_{10}^{16}(\al,\be),\; \mu_{10}^{26}(0,\be),\;  \mu_{10}^{27},\;
\mu_{10}^{29}(\al),\; \mu_{10}^{32}(-\frac{1}{2},\be),$ 
$\mu_{10}^{33}(-\frac{1}{2}),\;
\mu_{10}^{36}(0,\be),\;\mu_{10}^{39}(\al),\;
\mu_{10}^{40},\;\mu_{10}^{45}(\al),\;\mu_{10}^{46},\; \mu_{10}^{48}$ or

\begin{align*}
\mu_{10}^1(\al,\be) & \colon \quad (\al+2)(\al+1)(\al-1)p(1,\al) \neq 0 \\
\mu_{10}^2(\al,\be) & \colon \quad (\al+2)(\al+1)(\al-1)p(1,\al) \neq 0 \\
\mu_{10}^3(\al) & \colon \quad (\al+2)(\al+1)(\al-1)p(1,\al) \neq 0 \\
\mu_{10}^{12}(\al,\be) & \colon \quad \al\neq 0 \\ 
\mu_{10}^{17}(\al,\be,\ga) & \colon \quad (2\be+1)(3\al+3\ga-4\ga^2) \neq 0 \\ 
\mu_{10}^{17}(\al,\be,\ga) & \colon \quad 2\be+1=3\al+3\ga-4\ga^2= 0 \\
\end{align*}
\begin{align*}
\mu_{10}^{18}(\al,\be,\ga,\de) & \colon \quad (2\be+\de)(3\al+3\ga-4\ga^2)\neq 0 \\ 
\mu_{10}^{18}(\al,\be,\ga,\de) & \colon \quad 2\be+\de=3\al+3\ga-4\ga^2=\ga(\ga-6)= 0, \be \neq 0 \\
\mu_{10}^{18}(\al,\be,\ga,\de) & \colon \quad \al=42,\be=0,\ga=6,\de=0 \\ 
\mu_{10}^{19}(\al,\be) & \colon \quad 3\al+3\be-4\be^2 \neq 0\\ 
\mu_{10}^{20}(\al,\be,\ga) & \colon \quad (2\be+\ga)(6\al+1)\neq 0 \\
\mu_{10}^{20}(\al,\be,\ga) & \colon \quad 2\be+\ga=6\al+1= 0 \\
\mu_{10}^{21}(\al,\be,\ga) & \colon \quad \ga(2\al+\be)\neq 0 \\
\mu_{10}^{22}(\al,\be,\ga) & \colon \quad \ga(2\al+\be)\neq 0 \\
\mu_{10}^{24}(\al,\be,\ga) & \colon \quad 2\al+\be= 0, \al\neq 0 \\
\mu_{10}^{25}(\al,\be) & \colon \quad \al\neq 0 \\
\mu_{10}^{28}(\al,\be) & \colon \quad 2\be +1\neq 0 \\
\mu_{10}^{31}(\al,\be,\ga) & \colon \quad 2\al +\be= 0, \al\neq 0 \text{ or } \ga \neq 0 \\
\end{align*}
\end{cor}

\section{Symplectic filiform CNLAs of dimension $n\ge 12$}

For $n\ge 12$ there is no classification of symplectic filiform Lie algebras.
We will restrict ourselfs to certain families of filiform Lie algebras $\Lg$ of dimension $n\ge 12$. 
Consider the following conditions on $\Lg$: 
 \begin{itemize}
\item[(a)] $\Lg$ contains no one-codimensional subspace $U \supseteq \Lg^1$
such that $[U,\Lg^1]\subseteq \Lg^4$.
\item[(b)] $\Lg^{\frac{n-4}{2}}$
is abelian, if $n$ is even.
\item[(c)] $[\Lg^1,\Lg^1]\subseteq \Lg^6$.
\end{itemize}  

These properties are isomorphism invariants.

\begin{defi}
Let $\CA_n^1$ denote the set of $n$-dimensional filiform laws whose algebras
satisfy the properties $(a),(b),(c)$. Denote by $\CA_n^2$ the set of
$n$-dimensional filiform laws whose algebras satisfy $(a),(b)$, but {\it not}
$(c)$. Finally, for $n$ even, denote by $\CA_n^3$ the set of $n$-dimensional 
filiform laws whose algebras satisfy $(a)$ but not $(b)$. 
\end{defi}

The above properties of $\Lg$ can be expressed in terms of the 
corresponding structure constants $\al_{k,s}$. It is easy to verify the
following (use \eqref{lie}, $(5)$):

\begin{itemize}
\item[] $\al_{2,5}\ne 0$, if and only $\Lg$ satisfies property $(a)$.
\item[] $\al_{\frac{n}{2},n}=0$, if and only if $\Lg$ satisfies
property $(b)$.
\item[] $\al_{3,7}=0$, if and only if $\Lg$ satisfies property $(c)$.
\end{itemize}     

If $\Lg$ satisfies property $(a)$ we may change the adpated basis so that
it stays adapted and 
$$\al_{2,5}=1.$$
In fact, we may take $f\in GL(\Lg)$ defined by $f(e_1)=ae_1, f(e_2)=be_2$
and $f(e_i)=[f(e_1),f(e_{i-1})]$ for $3\le i\le n$ with suitable 
nonzero constants $a$ and $b$.\\

\begin{prop}
Suppose that $\Lg$ is a filiform Lie algebra of dimension $n\ge 12$ satisfying
properties $(a),(b)$. Hence we may assume for its law $\la\in \CA_n$ that
$\al_{2,5}=1$, and $\al_{\frac{n}{2},n}=0$ if $n$ is even. Then the Jacobi
identity implies that 

\begin{equation*}
(\al_{3,7},\al_{4,9},\al_{5,11})=\begin{cases}
(0,0,0) & \text{if $\la\in \CA_n^1$}\\[0.2cm]
(\frac{1}{10}, \frac{1}{70},\frac{1}{420}) & \text{if $\la\in \CA_n^2$}\\[0.2cm]
\end{cases}
\end{equation*}
\end{prop}

\begin{proof}
Let $(e_1,\ldots, e_n)$ be an adapted basis of $\Lg$, the Lie brackets
with respect to this basis being given by $\eqref{lie},(5)$. 
 Let $J(e_i,e_j,e_k)=0$ denote the Jacobi identity with $e_i,e_j,e_k$.
Let $J(i,j,k,l)$ be the coefficient of $e_l$ in $J(e_i,e_j,e_k)$.
If $n\ge 12$ then we have the conditions
$J(2,3,4,9)=J(2,4,5,11)=J(3,4,5,12)=0$ which are given by the following
equations:
\begin{gather*}
\al_{4,9}(2+\al_{3,7}) - 3\al_{3,7}^2=0\\
\al_{5,11}(2-\al_{3,7}-\al_{4,9}) + 2\al_{4,9}(3\al_{4,9}-2\al_{3,7})=0\\
3\al_{5,11}(\al_{3,7}+\al_{4,9}) -4\al_{4,9}^2=0
\end{gather*}
It is not difficult to see that there are precisely two solutions:
$(\al_{3,7},\al_{4,9},\al_{5,11})= (0,0,0)$ or
$(\al_{3,7},\al_{4,9},\al_{5,11}) = (\frac{1}{10},\frac{1}{70},
\frac{1}{420})$. 
Indeed, the first equation implies
$\al_{4,9}=3\al_{3,7}^2/(2+\al_{3,7})$. If we
substitute that into the other equations we obtain 
$\al_{3,7}(10\al_{3,7}-1)=0$.
\end{proof}     

\begin{prop}\label{twelve}
Let $\Lg$ be a filiform Lie algebra of dimension $12$ with law
$\la\in\CA_{12}$. Then \\
\begin{equation*}
H^2(\Lg,k)=\begin{cases}
\s\{[\om_1],[\om_2],[\om]\} & \text{if $\la\in \CA_{12}^1$}\\[0.2cm]
\s\{[\om_1],[\om_2],[\be]\} & \text{if $\la\in \CA_{12}^2$}\\[0.2cm]
\s\{[\om_1],[\om_2]\} & \text{if $\la\in \CA_{12}^3$}\\[0.2cm]
\end{cases}
\end{equation*}
If $\la\in \CA_{12}^2$ then $\Lg$ is symplectic since $\det (\be)\ne 0$.
In the other two cases $\Lg$ is not symplectic.
\end{prop}

\begin{proof}
Let $\la\in \CA_{12}^1$ and the Lie brackets of $\Lg$ being given by
$\eqref{lie},(5)$ with $21$ scalars $\al_{k,s}, (k,s)\in \CI_{12}$.
We have $\al_{2,5}=1, \, \al_{3,7}=\al_{6,12}=0$ and  
polynomial equations in the parameters $\al_{k,s}$
given by the Jacobi identity. However since we use 
an adapted basis, these equations are quite simple. The Jacobi identity
is satisfied if and only if\\
\begin{align*}
\al_{4,9} & =\al_{4,10}=\al_{5,11}=\al_{5,12}=0\\
\al_{4,11} & = 2\al_{3,8}^2 \\
\al_{4,12} & = -\frac{1}{2}\left[3\al_{4,11}(\al_{2,6}+\al_{3,8})-
9\al_{3,9}\al_{3,8}\right] 
\end{align*}       
Hence the parameters $\al_{2,6},\ldots,\al_{2,12}$ and $\al_{3,8},\ldots,
\al_{3,12}$ are arbitrary.\\
Now a standard computation yields the second scalar cohomology as above.
Here $\om$ is a $2$-cocycle with
 \begin{align*}
\om(e_1\wedge e_{12})& = 1\\
\om(e_2\wedge e_4)& = \al_{2,12}\\
\om(e_2\wedge e_6)& = \al_{2,10}-2\al_{3,12}\\
\vdots  \hspace{0.7cm}  & = \quad \vdots \\
\om(e_4\wedge e_7)& = 2\al_{3,8}^2\\
\end{align*}   
It is easy to see that all linear combinations of $\om_1,\om_2$ and $\om$
are degenerate. In fact, the vector $(0,\ldots,0,1,6\al_{3,8}-\al_{2,6},0)^t$
always belongs to the kernel of the representing matrix. 
Hence $\Lg$ is not symplectic.
A similar computation is done for the other two cases.
\end{proof}

\begin{cor}
If $\la\in \CA_{12}^2$ such that $200 \al_{3,8}-27\al_{2,6}\neq 0$ then $\Lg$ is a
symplectic CNLA. 
\end{cor}

\begin{rem}
Let $\la\in \CA_{12}^2$. Then $\Lg$ is always a CNLA, except for the case
\begin{align*} 
\al_{3,8} & = 27\al_{2,6}/200, \\
\al_{3,9} & = (4000\al_{2,7}+243\al_{2,6}^2)/28000 \\
\al_{3,10} & = (560000\al_{2,8}+100000\al_{2,6}\al_{2,7}-30213 \al_{2,6}^3)/3920000 \\
\al_{3,11} & = f(\al_{2,6},\ldots ,\al_{2,9}) \\ 
\al_{3,12} & = g(\al_{2,6},\ldots ,\al_{2,10}) 
\end{align*}
with certain polynomials $f,g \in \Q[\al_{2,6},\ldots ,\al_{2,10}]$.
\end{rem}

For $\la\in \CA_{n}^2, n\ge 13$ we have two different cases for the cohomology. 
Denote by $\CA_{n,1}^2$ the subset of laws satisfying
$\al_{3,n-4}=P_n(\al_{k,s})$, where $P_n$ is a certain polynomial with
rational coefficients in the variables $\al_{k,s}$ where 
$k=2, \,6\le s\le n$ and $k=3,\,n-4\le s \le n$. 
Denote by $\CA_{n,1}^2$ the subset of laws which do not satisfy this
polynomial equation. For $n=14$ the polynomial $P_{14}$ is given by:\\
\begin{equation*}
\begin{split}  
P_{14} & = (482832810500a_{3,8}^3 - 157196008500a_{3,8}^2a_{2,6} + 
2223828750a_{3,8}a_{2,7}\\
 & \hskip0.05 cm + 16180336845a_{3,8}a_{2,6}^2 + 186801615a_{2,8} - 
266859450a_{2,7}a_{2,6}\\
 & \hskip0.05 cm - 517476276a_{2,6}^3)/1307611305\\
\end{split}
\end{equation*} 

\begin{prop}\label{fourteen}
Let $\Lg$ be a filiform Lie algebra of dimension $14$ with law
$\la\in\CA_{14}$. Then \\
\begin{equation*}
H^2(\Lg,k)=\begin{cases}
\s\{[\om_1],[\om_2],[\om]\} & \text{if $\la\in \CA_{14}^1$}\\[0.2cm]
\s\{[\om_1],[\om_2],[\be]\} & \text{if $\la\in \CA_{14,1}^2$}\\[0.2cm]
\s\{[\om_1],[\om_2]\} & \text{if $\la\in \CA_{14,2}^2$}\\[0.2cm]
\end{cases}
\end{equation*}
If $\la\in \CA_{14,1}^2$ then $\Lg$ is symplectic since $\det(\be)\ne 0$. 
In the other two cases $\Lg$ is not symplectic.
\end{prop}

\begin{cor}
If $\la\in \CA_{14,1}^2$ such that $200 \al_{3,8}-27\al_{2,6}\neq 0$ then $\Lg$ is a
symplectic CNLA. 
\end{cor}

\begin{rem}
The result generalizes to higher dimensions. The cohomology has dimension
$2$ or $3$ and only the filiform algebras $\Lg$ with law
$\la\in \CA_{n,1}^2$ are symplectic. It is easy to see that
the algebras $\Lg$ with law $\la\in \CA_n^1, n\ge 12$ are not
symplectic, since $(0,\ldots,0,1,(n-6)\al_{3,8}-\al_{2,6},0)^t$
lies in the kernel of the matrix associated to every $\om\in Z^2(\Lg,k)$. 
The algebras with law $\la\in \CA_{n,1}^2$ are CLNAs except for the case
where $\al_{3,k}$ are given by certain polynomials in $\al_{2,l}$
with rational coefficients.
\end{rem}

\begin{rem}
The knowledge of $H^2(\Lg,k)$ can also be used to determine affine
filiform Lie algebras. If there exists a non-degenerate $\om\in Z^2(\Lg,k)$  
then $\Lg$ is symplectic, hence affine. However, it is enough to find
an affine class $[\om]\in H^2(\Lg,k)$ to ensure that $\Lg$ is affine, see \cite{BU1}. 
It is well known that all complex filiform Lie algebras of dimension  $n<10$ are affine.
In dimension $10$ however, there exist filiform algebras which are not affine. It turns out
that a law  $\la$ in  $\CA_{10,1}$ or $\CA_{10,4}$ is not affine, if it
belongs to a certain irreducible component of the variety of all nilpotent Lie algebra
laws of dimension $10$, such that the Lie algebra $\Lg/Z_2(\Lg)$ is characteristically nilpotent.
Here $Z_2(\Lg)$ denotes the second center of $\Lg$.
\end{rem}


\begin{thebibliography}{99}

\bibitem{BU1} D. Burde: {\it Affine cohomology classes for filiform Lie
algebras}. Contemp.\ Math.\ \textbf{262} (2000), 159--170. 

\bibitem{BU2} D. Burde: {\it Affine structures on nilmanifolds}. Int.\ J.\
of Math.\ \textbf{7} (1996), 599-616.  

\bibitem{BS} D. Burde, C. Steinhoff: {\it Classification of orbit closures
of $4$-dimensional complex Lie algebras}. J.\
of Algebra \textbf{214} (1999), 729--739. 

\bibitem{DAM} J.-M. Dardi\'e, A. M\'edina: {\it Alg\`ebres de Lie
kaehl\'eriennes et double extension}. J.\ Algebra \textbf{185} (1996), 
774--795.   

\bibitem{ELA} A. G. Ehlashvili: {\it Frobenius Lie algebras}. 
Funct. Anal. Appl. \textbf{16}, (1983), 326--328.

\bibitem{GJK} J.R. G\'omez, A. Jimenez-Merchan, Y. Khakimdjanov:
{\it Low-dimensional filiform Lie algebras}. J.\ Pure and Applied Algebra
\textbf{130} (1998), 133--158.   

\bibitem{GJKS} J.R. G\'omez, A. Jimenez-Merchan, Y. Khakimdjanov:
{\it Symplectic structures on filiform Lie algebras}. J.\ Pure and Applied Algebra
\textbf{156} (2001), 15--31. 

\bibitem{LAU} J. Lauret: {\it A distinguished compatible metric for geometric
structures on nilmanifolds}. Preprint (2004).

\bibitem{MIL1} D. V. Millionschikov: {\it Graded filiform Lie algebras and symplectic nilmanifolds}. 
Geometry, topology, and mathematical physics, Amer. Math. Soc. Transl. Ser. 2, \textbf{212} (2003),
259--279.

\bibitem{MIL2} D. V. Millionschikov: {\it Deformations of graded nilpotent Lie algebras
and symplectic structures}. Preprint arXiv. \textbf{0305057}, (2003).

\bibitem{OOM} A. I. Ooms: {\it On Frobenius Lie algebras}. Comm.\
Alg.\ \textbf{8} (1980), 13--52.

\bibitem{PAR} S. E. Parkhomenko: {\it Quasi-Frobenius Lie algebras construction of $N=4$ superconformal 
field theories}. Mod. \ Phys. \ Lett. A \textbf{11} (1996) , No.6, 445--461..

\bibitem{STO} A. Stolin: {\it Rational solutions of the classical
Yang-Baxter equation and quasi Frobenius Lie algebras}.
J.\ Pure Appl.\ Algebra \textbf{137} (1999), 285--293.   

\bibitem{MIL} J. Milnor: {\it On fundamental groups of complete affinely
flat manifolds}. Advances in Math.\ \textbf{25} (1977), 178--187.


\end{thebibliography}
\end{document}